\documentclass[10pt,a4paper,twocolumn]{article}

\usepackage{amsmath,amsfonts,amssymb,dsfont,graphicx,color}
\newcommand{\rea}{\mathbb{R}}
\newcommand{\nat}{\mathbb{N}}

\newtheorem{theorem}{\textbf{Theorem}}

\newtheorem{propo}{\textbf{Proposition}}

\begin{document}
\title{Consensus on the average in arbitrary directed network topologies with time-delays}

\author{Mehran Zareh$^*$, Carla Seatzu$^*$, Mauro Franceschelli$^*$
	\\~\\
	$^{*}$ Department of Electrical and Electronic Engineering, University of Cagliari, Italy \\
	\{mehran.zareh,mauro.franceschelli,seatzu\}@diee.unica.it}

\thanks{The research leading to these results has received funding from the European Union Seventh Framework Programme [FP7/2007-2013]  under grant agreement n 257462 HYCON2 Network of excellence.}

\maketitle

\begin{abstract}
In this preliminary paper we study the stability property of a consensus on the average algorithm in arbitrary directed graphs with respect to communication/sensing time-delays. The proposed algorithm adds a storage variable to the agents' states so that the information about the average of the states is preserved despite the algorithm iterations are performed in an arbitrary strongly connected directed graph. We prove that for any network topology and choice of design parameters the consensus on the average algorithm is stable for sufficiently small delays. We provide simulations and numerical results to estimate the maximum delay allowed by an arbitrary unbalanced directed network topology.
\end{abstract}

\section{Introduction}

The consensus problem in multi-agent systems consists in the design of a coupling law between dynamical systems (agents) such that the state of each one converges to the same value in absence of external reference signals. Multi-agent systems are considered to be \emph{complex} systems since the pattern of interconnections between agents is often arbitrary and unknown at the controller design stage. This clearly makes challenging the design of interaction rules between agents that exploit only local information. For these reasons agents modeled by simple single integrators or second order systems are usually investigated. One of the major works from which we take inspiration is the one by \cite{olfati2004consensus} where the consensus problem for networks of first order agents for switching topologies or time-delays is investigated. In this paper the authors prove that simple averaging local interaction rules can achieve consensus on the average, i.e., the state of each agent converges to the average of the initial states only if the directed graph that encodes the network topology is strongly connected and balanced (each agent receives and sends information to the same number of agents). They also explored the consensus problem in the case of time-delays for undirected network topologies.

Since then several authors have explored ways to design consensus on the average algorithms that work on general directed graphs not necessarily balanced. In \cite{FEG:2008,FGS:2009} the idea to use an augment state space to add robustness to a networked system represented by an undirected graph that executes a consensus algorithm was proposed. The proposed algorithms aim at recovering the correct network average once malicious or faulty agents have been removed from the network.

In \cite{FGS:2009,fra2010} a discrete time consensus on the average algorithm for arbitrary strongly connected directed graph based on asynchronous state updates was presented, based on the idea to augment the state of each agent with an additional variable to preserve the information about the initial average of the states in the network. Simulations were used to characterize the convergence properties and the performance of the algorithm.

In \cite{cai2012average} a discrete time consensus on the average algorithm based on additional state variables was characterized in terms of a tuning parameter. It was proven that there always exist sufficiently small values of such tuning parameter so that the proposed algorithm converges to the average of the initial state in arbitrary strongly connected directed graphs.

In \cite{caiquantized2011} a consensus on the average algorithm for agents with uniformly quantized state and asynchronous updates in directed graphs was proposed.

In \cite{hadj2012} the authors address the control of distributed energy resources by developing a consensus on the average protocol based on the so called \emph{ratio} consensus. Their algorithm is based on two independent distributed dynamical systems, one with arbitrary initial conditions and one with predetermined initial values. The authors consider time-varying network topologies described by directed graphs and show that for each agent the ratio of the output of these two dynamical systems converges to the average of the initial states.

In \cite{Yin2010,Yin2011} the \emph{Corrective Consensus} algorithm is proposed. It consists in a local state update rule where each agent keeps track of several additional variables corresponding to the number of its neighbors which are used to periodically steer the average of the network state to the correct value corresponding to the average of the states at the initial instant of time.

In \cite{Scaglione2009} the broadcast gossip algorithm is proposed. This algorithm is based upon discrete time and asynchronous state updates with directed information flow, it makes each agent agree upon a random variable whose expectation is the average of the initial states.

Most of the literature on consensus on the average in directed graphs deals with methods and techniques to achieve consensus on the average in networks of agents described by single integrators. On the other hand the literature on consensus with time-delays in directed graphs usually deals with the problem of making the state of each agent converge to the same value which can be time-varying and not related to the initial state of the network in an explicit way.

In \cite{yu2010some} necessary and sufficient conditions for convergence of second-order multi-agent systems with velocity feedback are given and the effect of time-delays in directed graphs is characterized. In \cite{yu2010some} the consensus value is arbitrary.

In \cite{sun2009consensus} several instances of consensus problems with time-delays are investigated. In particular the cases of switching directed topologies, packet data dropouts, and finite time consensus are all characterized separately by considering the effect of time-delays for the achievement of consensus on an arbitrary value.

In this preliminary paper we propose a continuous time consensus algorithm inspired from the discrete time algorithms in \cite{FGS:2009,fra2010,cai2012average}.
We consider a description in continuous time to describe a network of $n$ vehicles  with a local interaction rule that controls the instantaneous speed of each vehicle.
Then we extend the proof method in \cite{cai2012average} to the case at hand and study the convergence properties of the resulting system considering a time-delay in the state update of each agent. We finally provide simulation results to corroborate the theoretical analysis.

The main contributions of this paper can be summarized in the following three items.
\begin{itemize}
\item We provide a continuous time version of a consensus on the average algorithm for arbitrary directed strongly connected graphs derived from results in \cite{FGS:2009,fra2010} and \cite{cai2012average}.
\item We provide a characterization of the convergence properties of the algorithm with respect to time-delays.
\item We present simulations to characterize numerically the performance of the proposed protocol with respect to different time-delays and tuning parameters.
\end{itemize}

The paper is structured as follows. In Section~\ref{section:Prob_for} we introduce the notation and recall some preliminaries of algebraic graph theory. In Section~\ref{section:protocol} we introduce a consensus on the average protocol and the corresponding model considering time-delays. In Section~\ref{section:convergence} we characterize the convergence properties of the proposed algorithm with respect to time-delays. In Section~\ref{section:example} we corroborate the theoretical analysis with a numerical example and simulations. Concluding remarks are finally given in Section~\ref{section:conclusion}.

\section{Preliminaries}\label{section:Prob_for}
We now introduce some basic concepts of algebraic graph theory. The pattern of interactions among agents (nodes) is represented by a directed graph $\mathcal{G}=(\mathcal{V}, \mathcal{E})$ where $\mathcal{V}=\{1,\ldots,n\}$ is the set of nodes and $\mathcal{E}\subseteq\{\mathcal{V}\times \mathcal{V}\}$ is the set of directed edges. Node $i$ receives information from node $j$ only if $(i,j)\in \mathcal{E}$.

The sets of in-neighbors and out-neighbors which represent the set of nodes from which each node receives or sends information are defined as $\mathcal{N}_{i,in}=\{ j \in \mathcal{V},(i,j)\in \mathcal{E} \} $ and $\mathcal{N}_{i,out}=\{j \in \mathcal{V},(j,i)\in \mathcal{E}\} $, respectively.

The topology of graph $\mathcal{G}$ is encoded by the so-called adjacency matrix, an $n \times n$ matrix $A=[a_{ij}]$ whose elements are $a_{ij}=1$ if $(i,j)\in \mathcal{E}$ and $0$ otherwise.

Let the in-degree and out-degree of node $i$ be respectively $\delta_{i,in}=|\mathcal{N}_{i,in}|$ and $\delta_{i,out}=|\mathcal{N}_{i,out}|$, where by $|\mathcal{N}_{i,in}| \ \left(|\mathcal{N}_{i,out}|\right)$ we denote the cardinality of set $\mathcal{N}_{i,in} \ \left(|\mathcal{N}_{i,out}|\right)$.
Let $D_{in}=\emph{diag}(\delta_{1,in}, \ldots,\delta_{n,in})$ and $D_{out}=\emph{diag}(\delta_{1,out}, \ldots,\delta_{n,out})$ be $n \times n$ diagonal matrices whose diagonal elements are respectively the in-degree and out-degree of the nodes.

We now define the in-Laplacian matrix as $\mathcal{L}_{in}=D_{in}-A$ and the out-Laplacian matrix as $\mathcal{L}_{out}=D_{out}-A$. The in-Laplacian and out-Laplacian matrices of a directed graph have several structural properties. Due to the Gershgorin Circle Theorem applied to the rows of the in-Laplacian or the columns of the out-Laplacian it is possible to show that both matrices have eigenvalues with non-negative real part for any graph $\mathcal{G}$. By construction matrices $\mathcal{L}_{in}$ and $\mathcal{L}_{out}$ have at least one null eigenvalue because either the row sum or the column sum is zero. Furthermore, let $\mathbf{1}_n$ and $\mathbf{0}_n$ be respectively the $n$-elements vectors of ones and zeros, then $\mathcal{L}_{in}\mathbf{1}=\mathbf{0}$ and $\mathbf{1}^T\mathcal{L}_{out}=\mathbf{0}^T$. If $\mathcal{G}$ is strongly connected, i.e., there exists a directed path that connects any pair of nodes in $\mathcal{V}$, then the algebraic multiplicity of the null eigenvalue of both $\mathcal{L}_{in}$ and $\mathcal{L}_{out}$ is one.
%
%
%
\section{Consensus on the average protocol}\label{section:protocol}

We now introduce a consensus protocol stated in continuous time that takes inspiration from protocols appeared in \cite{fra2010} and \cite{cai2011average} in a discrete time setting.  In the protocol under consideration each agent is a single integrator with an additional state variable called \emph{surplus} or \emph{storage}. This additional variable is used to preserve information about the average value of the agents' states at the initial instant of time, that is a time-varying quantity in directed graphs that are not balanced, i.e., graphs in which the in-degree and out-degree of each node are not necessarily equal.

The local state update rule implemented by each node is the following:
\begin{equation} \label{consensusprotocol}
\left\{
\begin{array} {ll}
\dot{x}_i(t)=&-\sum_{j \in \mathcal{N}_{i,in}}\left(x_i(t)-x_j(t)\right)+ \varepsilon z_i(t),\\
\dot{z}_i(t)=&\sum_{j \in \mathcal{N}_{i,in}}\left(x_i(t)-x_j(t)\right) \\
 & -\sum_{j \in \mathcal{N}_{i,in}}\left(z_i(t)-z_j(t)\right) \\
 &-\left(\varepsilon -\delta_{i,in}+\delta_{i,out}\right) z_i(t),
\end{array}\right.
\end{equation}
where $x_i,z_i \in \rea$ are the states of agent $i$ and $\varepsilon\in \rea^+$ is a tuning parameter of the algorithm. It is clear that to implement protocol~\eqref{consensusprotocol} each agent requires only relative state information with respect to variable $x_i$, absolute state information with respect to variable $z_i$, and knowledge of its own out-degree.

The network dynamics that emerges when each agent implements the local state update rule in eq.~\eqref{consensusprotocol} can be formulated in matrix form as follows:
\begin{equation} \label{withoutdelay}
\left[ \begin{array}{c}
 \dot{x}(t) \\
\dot{z}(t)
\end{array} \right]=\left[ \begin{array}{cc}
-\mathcal{L}_{in} & \varepsilon I \\
 \mathcal{L}_{in} & -\mathcal{L}_{out}-\varepsilon I
\end{array} \right] \left[ \begin{array}{c}
 {x}(t) \\
{z}(t)
\end{array} \right]
\end{equation}

where $x=\left[x_1, x_2,\ldots,x_n\right]$ and  $z=\left[z_1, z_2, \ldots, z_n\right]$ are a compact representation of the agents' state.

{ The proposed local interaction scheme can be interpreted as a network of $n$ vehicles each modeled as a continuous time single integrator $\dot{x}_i=u_i(t)$ where each $x_i(t)$ represents a position in space and  variables $z_i(t)$ are software variables which enable the interaction scheme to converge to the initial average position.}



In this paper we study protocol~\eqref{withoutdelay} under the assumption that communication/sensing delays affect the multi-agent system. The network dynamics are thus described by

\begin{equation} \label{withdelay}
\begin{array}{c}
 \left[ \begin{array}{c}
  \dot{x}(t) \\
 \dot{z}(t)
 \end{array} \right]= M(\varepsilon) \left[ \begin{array}{c}
  {x}(t-\tau) \\
 {z}(t-\tau)
 \end{array} \right]
\end{array}
\end{equation}

with $x(\theta)=x_0, \quad z(\theta)=z_0, \quad \ -\tau \le \theta \le 0,$ where

\begin{equation}\label{M}
M(\varepsilon)=\left[ \begin{array}{cc}
-\mathcal{L}_{in} & \varepsilon I \\
 \mathcal{L}_{in} & -\mathcal{L}_{out}-\varepsilon I
\end{array} \right]
\end{equation}
and $\tau \in \rea^+$ denotes a time-delay. We study system~\eqref{withdelay} in the approximation that the delay for all the agent is the same.

\section{Convergence properties}\label{section:convergence}

In thiconss section we study the convergence properties of system~\eqref{withdelay}.

We preliminary observe that by construction matrix $M(\varepsilon)$ satisfies $\left[\mathbf{1}_n^T \ \mathbf{1}_n^T\right]M(\varepsilon)=\left[\mathbf{0}_n^T \ \mathbf{0}_n^T\right]$ for any $\varepsilon \in \rea$. Therefore, since $$\mathbf{1}_n^T\dot{x}(t)+\mathbf{1}_n^T\dot{z}(t)=0, \quad \forall t\geq 0$$ it holds
\begin{equation}\label{eq1}
\mathbf{1}_n^Tx(t)+\mathbf{1}_n^Tz(t)=\mathbf{1}_n^Tx(0)+\mathbf{1}_n^Tz(0), \quad \forall t\geq 0.
\end{equation}

Now consider matrix $M(\varepsilon)$ for $\varepsilon=0$, namely
\begin{equation}\label{mdizero}
M(0)=\left[\begin{array}{cc}
 -\mathcal{L}_{in}& 0  \\
 \mathcal{L}_{in}& -\mathcal{L}_{out}
\end{array} \right].
\end{equation}
It is clear that since matrix $M(0)$ is a $2n \times 2n$ block lower triangular matrix it has $2n$ eigenvalues equal to the eigenvalues of matrices $ - \mathcal{L}_{in}$ and $ -\mathcal{L}_{out}$. If graph $\mathcal{G}$ is strongly connected, then $M(0)$ has one null eigenvalue with algebraic multiplicity $2$ and geometric multiplicity $2$, all other eigenvalues have strictly negative real part.

{ In the following we denote as $\lambda_i(0)$, $i=1,\ldots,2n$, the eigenvalues of matrix $M(0)$ and assume that $$0=\lambda_1(0)=\lambda_2(0)> \Re(\lambda_3)(0) \geq \ldots \geq \Re(\lambda_{2n}(0)).$$ Eigenvalues of matrix $M(\varepsilon)$ are denoted as $\lambda_i(\varepsilon)$, $i=1,\ldots,2n,$ and ordered as $\Re(\lambda_1)(\varepsilon)\geq  \ldots \geq \Re(\lambda_{2n}(\varepsilon)).$
}

We now prove some properties of the eigenvalues of matrix $M(\varepsilon)$ for small values of $\varepsilon>0$, that can be derived from the results in \cite{cai2011average}.

\begin{propo}\label{propointro}
Let matrix $M(\varepsilon)$ be defined as in eq.~\eqref{M}. If $\mathcal{G}$ is strongly connected, there exists $\bar{\varepsilon}\in \rea^+$ such that if $\varepsilon \in \left(0,\bar{\varepsilon}\right]$ then $M(\varepsilon)$ has one null eigenvalue and $2n-1$ eigenvalues with strictly negative real part.

\textit{Proof}:  Matrix $M(\varepsilon)$ depends smoothly on parameter $\varepsilon \geq 0$, therefore if eigenvalues $\lambda_3(0),\ldots,\lambda_{2n}(0)$ of $M(0)$ have strictly negative real part, there exists $\bar{\varepsilon}>0$ such that if $\varepsilon\in \left[0,\bar{\varepsilon}\right]$ then for $i=3,\ldots,2n,$ it holds $\Re(\lambda_i(\varepsilon))<0$. Therefore, as in \cite{cai2011average}, we only have to show that for $\varepsilon$ sufficiently small { it is $\lambda_1(\varepsilon)=0$ and} $\Re(\lambda_2(\varepsilon))<0$.

Since the null eigenvalue of $M(0)$ is semi-simple\footnote{An eigenvalue is semi-simple if its algebraic and geometric multiplicity are equal.} and $Rank(M(0))=2n-2$ it has two linearly independent right eigenvectors $r_1, r_2$ and left eigenvectors $l_1, l_2$ corresponding to the null eigenvalue.
It holds
\begin{equation}\label{M_2parts}
M'=\frac{d M(\varepsilon)}{d\varepsilon}= \left[\begin{array}{cc}
 0& I  \\
 0& -I
\end{array}\right].
\end{equation}

Then, as shown in \cite{cai2011average},  $d\lambda_1(\varepsilon)/d\varepsilon|_{\varepsilon=0}$ and $d\lambda_2(\varepsilon)/d\varepsilon|_{\varepsilon=0}$ are the eigenvalues of the following matrix

\begin{equation} \label{mattemp}
\left[\begin{array}{cc}
 l_1^T M' r_1 & l_1^T M' r_2  \\
l_2^T M' r_1 & l_2^T M' r_2
\end{array} \right].
\end{equation}
If graph $\mathcal{G}$ is strongly connected then $l_1=\alpha_1 \mathbf{1}_{2n}$ and $r_1=\alpha_2 \left[\mathbf{1}_n^T,\mathbf{0}_{n}^T\right]$ where $\alpha_1,\alpha_2\in \rea$ can be chosen such that $l_1^Tr_1=1$. By substituting $l_1$ and $r_1$ in \eqref{mattemp} it can be shown by simple computations that

$$d\lambda_1(\varepsilon)/d\varepsilon|_{\varepsilon=0}=0, \quad d\lambda_2(\varepsilon)/d\varepsilon|_{\varepsilon=0}=l_2^T M' r_2.$$

{ The first equality enables us to conclude that for sufficiently small values of $\varepsilon$, it is $\lambda_1(\varepsilon)=0$. }

Now, let $\nu_{r,out}$ be the right eigenvector corresponding to the null eigenvalue of matrix $\mathcal{L}_{out}$ and $\nu_{l,in}$ be the left eigenvector corresponding to the null eigenvalue of matrix $\mathcal{L}_{in}$. It is possible to verify by substitution that we can choose $r_2=\left[\mathbf{0}_n^T,\nu_{r,out}^T\right]$ and $l_2=\left[\nu_{l,in}^T,\mathbf{0}_n^T\right]$. Therefore,

$$d\lambda_2(\varepsilon)/d\varepsilon|_{\varepsilon=0}=-\nu_{l,in}^T\nu_{r,out}.$$

Since $\mathcal{L}_{in}$ and $\mathcal{L}_{out}$ are Metzler matrices (\cite{berman1979nonnegative}), the eigenvectors $\nu_{l,in}$ and $\nu_{r,out}$ corresponding to the null eigenvalue have only positive elements. Therefore
$$d\lambda_2(\varepsilon)/d\varepsilon|_{\varepsilon=0}=-\nu_{l,in}^T\nu_{r,out}<0$$
and $\lambda_2(\varepsilon)<0$ for $\varepsilon>0$ sufficiently small, thus proving the statement. \hfill $\square$

\end{propo}

We are now ready to study the stability of system~\eqref{withdelay} with respect to time-delays. Let $Y(s)=\left[X(s)^T \ Z(s)^T\right]^T$ denote the Laplace transform of $y(t)=\left[x(t)^T \ z(t)^T\right]^T$. Then the Laplace transform of system~\eqref{withdelay} is
$$Y(s)=\left(sI-M(\varepsilon)e^{-s\tau}\right)^{-1}Y(0)$$
and the stability property of system~\eqref{withdelay} depends upon the roots of the quasi-polynomial
\begin{equation}\label{charpoly}
det\left(sI-M(\varepsilon)e^{-s\tau}\right).
\end{equation}
By simple manipulations it holds
\begin{equation}\label{eqr1}
det\left(sI-M(\varepsilon)e^{-s\tau}\right)=e^{-2ns\tau}det\left(se^{s\tau}I-M(\varepsilon)\right)
\end{equation}
thus the roots of \eqref{charpoly} correspond to the solutions of
\begin{equation}\label{eq123}
se^{s \tau}=\lambda_i(\varepsilon), \quad i=1,\ldots,2n.
\end{equation}

\begin{theorem} \label{teorema1}
Let matrix $M(\varepsilon)$ be defined as in eq.~\eqref{M} and $\varepsilon \in \left(0,\bar{\varepsilon}\right]$ as in Proposition~\ref{propointro}. If $\mathcal{G}$ is strongly connected and
\begin{equation}\label{eq13}
\tau\leq \tau_c(\varepsilon)=\min_{i=2,\ldots,2n} \frac{\theta_i(\varepsilon)-\frac{\pi}{2}}{R_i(\varepsilon)},
\end{equation}
where $R_i(\varepsilon)=|\lambda_i(\varepsilon)|$ and $\theta_i(\varepsilon)=\angle \lambda_i(\varepsilon)$ with $\lambda_i(\varepsilon)$ the $i$-th eigenvalue of $M(\varepsilon)$, then the roots of
\begin{equation}\label{eq41}
det\left(sI-M(\varepsilon)e^{-s\tau}\right)
\end{equation}
have all strictly negative real part except one in $s=0$.

\textit{Proof}: By Proposition~\ref{propointro} since $\mathcal{G}$ is strongly connected by assumption, there exists $\bar{\varepsilon}$ such that for $\varepsilon \in \left(0, \bar{\varepsilon}\right]$, $M(\varepsilon)$ has a single null eigenvalue and $2n-1$ eigenvalues with strictly negative real part. Since the roots of  eq.~\eqref{eq41} depend continuously on $\tau$ and for $\tau=0$ they coincide with the roots of $M(\varepsilon)$, we compute the smallest positive value of $\tau$, denoted as $\tau_c$, for which at least one non-null root crosses the imaginary axis. By eq.~\eqref{eq123}, assuming $s=j\omega$ it holds
$$j\omega e^{j\omega \tau}=R_i(\varepsilon)e^{j\theta_i(\varepsilon)}.$$
By simple manipulations the above equation can be rewritten as
$$j\omega=R_i(\varepsilon) \ cos(\theta_i(\varepsilon)-\omega \tau)+jR_i(\varepsilon) \ sin(\theta_i(\varepsilon) -\omega \tau),$$
therefore
$$\left\{\begin{array} {l}
R_i(\varepsilon) \ cos(\theta_i(\varepsilon)-\omega\tau)=0, \\
\omega=R_i(\varepsilon) \ sin(\theta_i(\varepsilon)-\omega \tau).
\end{array}
\right.$$
This implies that
$$\left\{\begin{array} {l}
\theta_i(\varepsilon)-\omega\tau=\displaystyle \frac{\pi}{2}+k\pi, \quad k\in \nat  \\
\omega=R_i(\varepsilon)\ sin(\displaystyle \frac{\pi}{2}+k\pi)=R_i(\varepsilon) (-1)^k.
\end{array}
\right.$$
Finally, considering only the top-half of the Gauss plane, $\theta_i(\varepsilon) \in \left(\frac{\pi}{2},\pi\right]$ for $i=1,\ldots,2n$. Thus

\begin{equation} \label{tau_c}
\begin{array}{ll}
\tau_c(\varepsilon) & = \displaystyle \min_{i=2,\ldots,2n} \displaystyle \min_{k\in \nat} \displaystyle \frac{\theta_i(\varepsilon)-\frac{\pi}{2}-k\pi}{R_i(\varepsilon)(-1)^k} \\~\\ & =\displaystyle \min_{i=2,\ldots,2n}  \displaystyle \frac{\theta_i(\varepsilon)-\frac{\pi}{2}}{R_i(\varepsilon)},
\end{array}
\end{equation}
proving the statement. \hfill $\square$
\end{theorem}

Next we give bounds on the maximum length of the time delay that ensures stability as function of known network parameters computed for $\varepsilon=0$. If the actual time delay is smaller than the proposed bound then we are sure that there exist $\varepsilon>0$ sufficiently small such that the system is stable and achieves consensus.

\begin{theorem}\label{coro}
Consider a multi-agent system that implements protocol \eqref{consensusprotocol} in graph $\mathcal{G}=\{\mathcal{V},\mathcal{E}\}$, with tuning parameter $\varepsilon>0$, initial condition $z(0)=\mathbf{0}_n$ and time-delay $\tau>0$. If $\mathcal{G}$ is strongly connected, there exists $\tilde{\varepsilon}$ such that if $\varepsilon \in \left(0,\tilde{\varepsilon}\right]$ and
$$\tau< \tilde{\tau} = \frac{1}{2 \bar{\delta}}\arctan \left(\frac{\Re \{\lambda_3(0)\}}{\bar{\delta}}\right)$$
where $$\bar{\delta}=\max_{i\in \mathcal{V}} \{\delta_{i,in},\delta_{i,out}\}$$ and $\lambda_3(0)$ is the rightmost non-null eigenvalue of matrix $M(0)$, then
$$\lim_{t \rightarrow \infty} x(t)=\frac{\mathbf{1}_n^Tx(0)}{n}\mathbf{1}_n.$$

\textit{Proof}:  By definition it holds
$$M(\varepsilon)=M(0)+\varepsilon M',$$
where $M'$ is defined as in eq.~\eqref{M_2parts}. Since $M(\varepsilon)$ can be seen as a perturbation of matrix $M(0)$ its eigenvalues depend continuously on parameter $\varepsilon$. This implies that the ratio in eq.~\eqref{eq13} can be bounded for an arbitrary small $\varepsilon$ as a function of the eigenvalues of $M(0)$. In particular, for $\varepsilon=0$ by the Gershgorin disc theorem applied to matrices $\mathcal{L}_{in}$ and $\mathcal{L}_{out}$ we have $R_i(\varepsilon)\leq \max_{i=1,\ldots,2n} |\lambda_i(\varepsilon)|\leq 2\max_{i\in \mathcal{V}} \{\delta_{i,in},\delta_{i,out}\}=2\bar{\delta}$, thus it holds
$$\begin{array}{ll}
\displaystyle \min_{i=2,\ldots,2n} \displaystyle \frac{\theta_i(\varepsilon)-\frac{\pi}{2}}{R_i(\varepsilon)} & \geq \displaystyle \frac{\min_{i=2,\ldots,2n} \theta_i(\varepsilon)-\frac{\pi}{2}}{\max_{i=2,\ldots,2n} R_i(\varepsilon)} \\~\\
                                                        & \geq \displaystyle \frac{1}{2 \bar{\delta}}\arctan \left(\min_{i=1,\ldots,2n} \frac{\Re (\lambda_i(\varepsilon))}{\Im (\lambda_i(\varepsilon))}\right).
\end{array}$$
Finally, since for $\varepsilon=0$, it is $\Im (\lambda_2(\varepsilon))=0$ and $\max_{i=1,\ldots,2n} |\Im (\lambda_i(\varepsilon))|\leq \bar{\delta}$, it holds
$$\min_{i=2,\ldots,2n} \frac{\theta_i(\varepsilon)-\frac{\pi}{2}}{R_i(\varepsilon)}\geq \frac{1}{2 \bar{\delta}}\arctan \left(\frac{\Re \{\lambda_3(\varepsilon)\}}{\bar{\delta}}\right).$$
Therefore, since by Theorem~\ref{teorema1} we may conclude that for $\tau\leq \tau_c(\varepsilon)$ all the roots of eq.~\eqref{charpoly} have strictly negative real part except one, this also holds for a sufficiently small value of $\varepsilon$ provided that
$$\tau<\frac{1}{2 \bar{\delta}}\arctan \left(\frac{\Re \{\lambda_3(0)\}}{\bar{\delta}}\right)=\tilde{\tau}\leq\tau_c(\varepsilon).$$

Therefore, the solutions $x(t)$ and $z(t)$ of system \eqref{withdelay} converge to the null space of matrix $M(\varepsilon)$, i.e.,
$$\lim_{t\rightarrow \infty}
\left[\begin{array} {c}
x(t)\\
z(t)
\end{array}\right]=r_1=\alpha\left[\begin{array} {c}
\mathbf{1}_n\\
\mathbf{0}_n
\end{array}\right].$$
Since $\mathbf{1}_n^Tx(t)+\mathbf{1}_n^Tz(t)=\mathbf{1}_n^Tx(0)+\mathbf{1}_n^Tz(0)$ for any $t\geq0$ we have that $$\alpha=\displaystyle \frac{\mathbf{1}_n^Tx(0)+\mathbf{1}_n^Tz(0)}{n}.$$ Since by assumption $z(0)=\mathbf{0}_n$, it holds
$$\lim_{t\rightarrow \infty} x(t)=\frac{\mathbf{1}_n^Tx(0)}{n}\mathbf{1}_n,$$
thus proving the statement. \hfill $\square$
\end{theorem}

\begin{figure}
\centering
\includegraphics[width=0.5\linewidth]{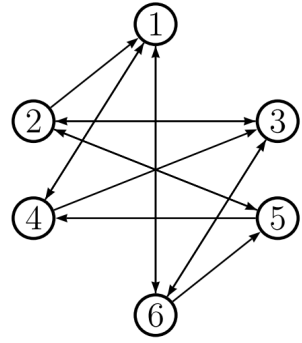}
\caption[Digraph]{The directed graph considered in Section~\ref{section:example}.}
\label{topology}
\end{figure}

\begin{figure}
\centering
\includegraphics[width=0.9\linewidth]{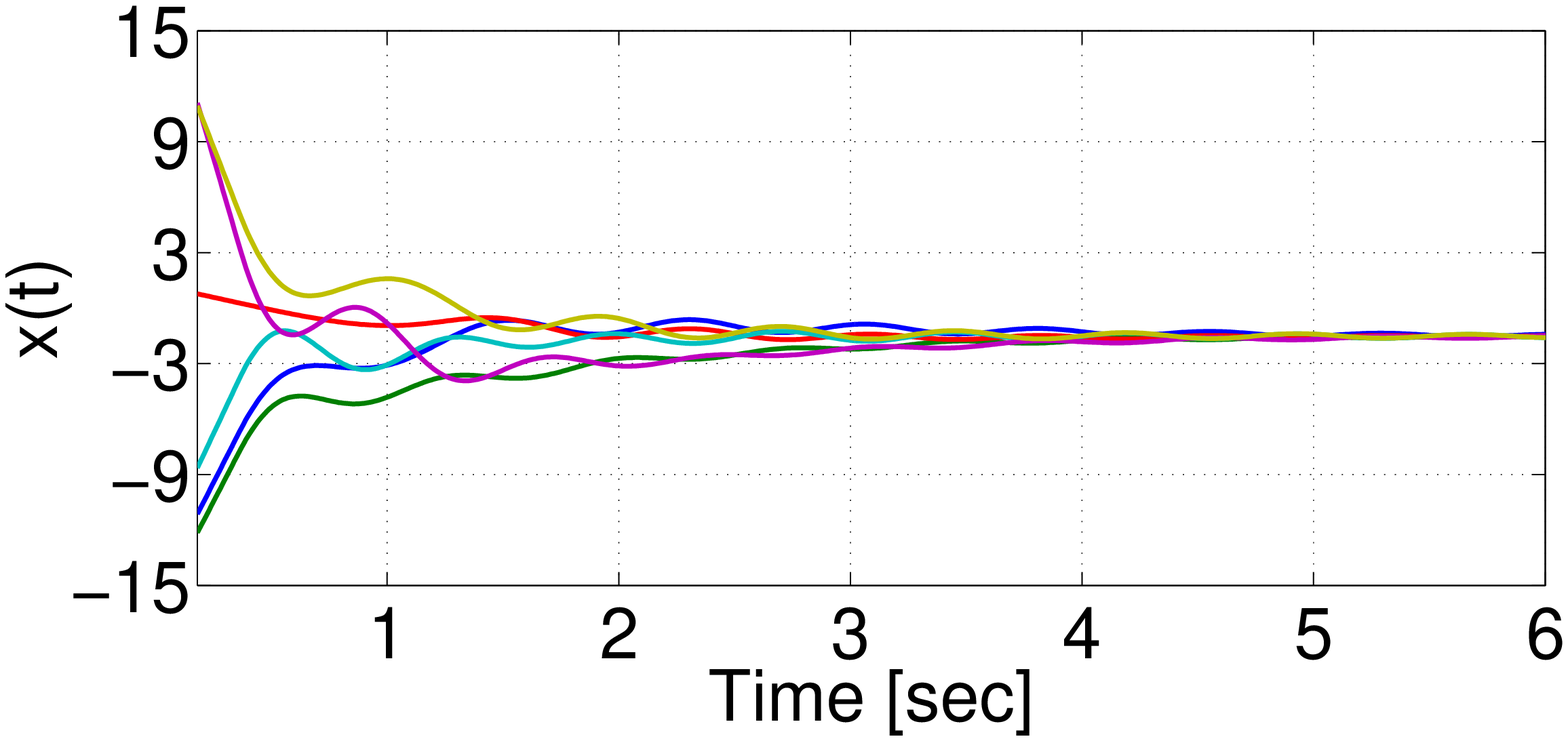}
\caption{Evolution of $x(t)$ for $\varepsilon=1.3$ and $\tau=0.19$.}\label{fig:examplex}
\end{figure}

\begin{figure}
\centering
\includegraphics[width=0.9\linewidth]{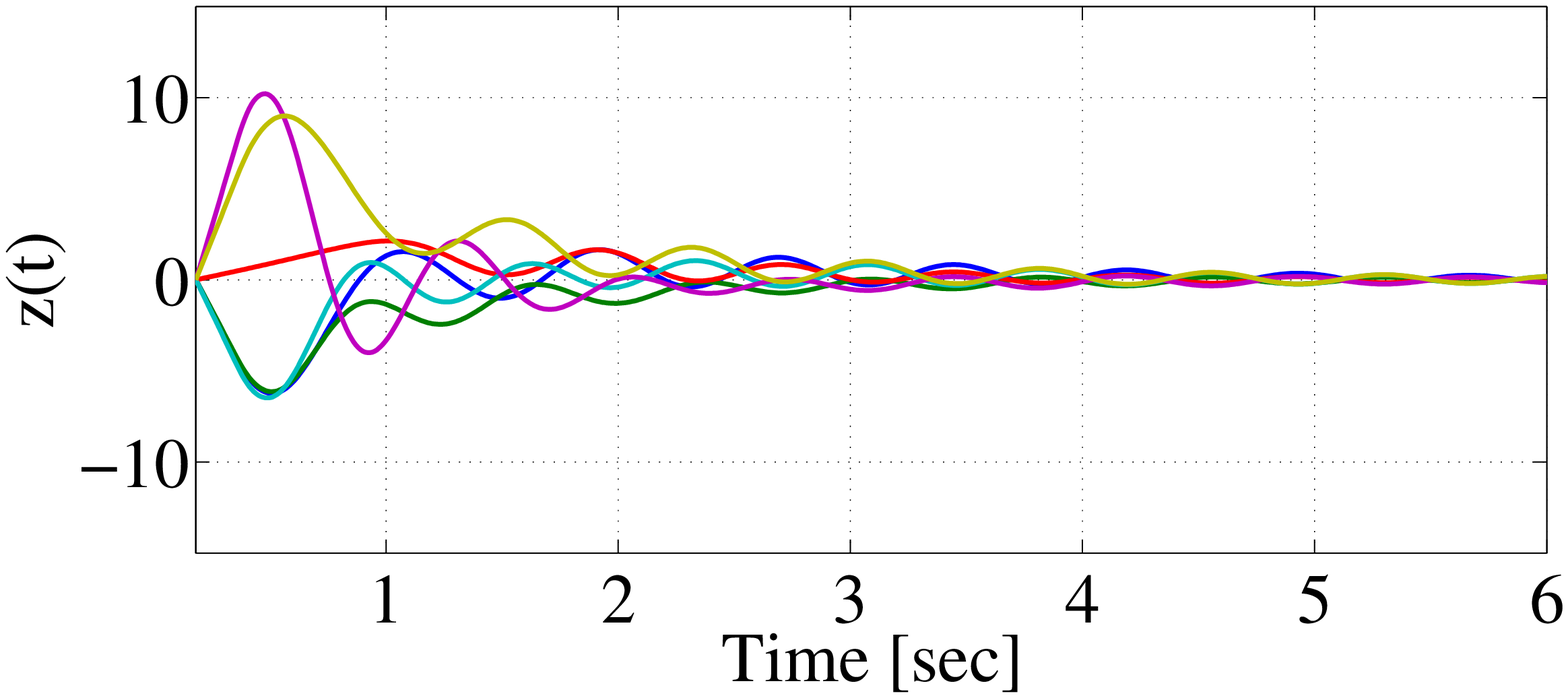}
\caption{Evolution of $z(t)$ for $\varepsilon=1.3$ and $\tau=0.19$.}\label{fig:examplez}
\end{figure}

\begin{figure}
\centering
\includegraphics[width=0.9\linewidth]{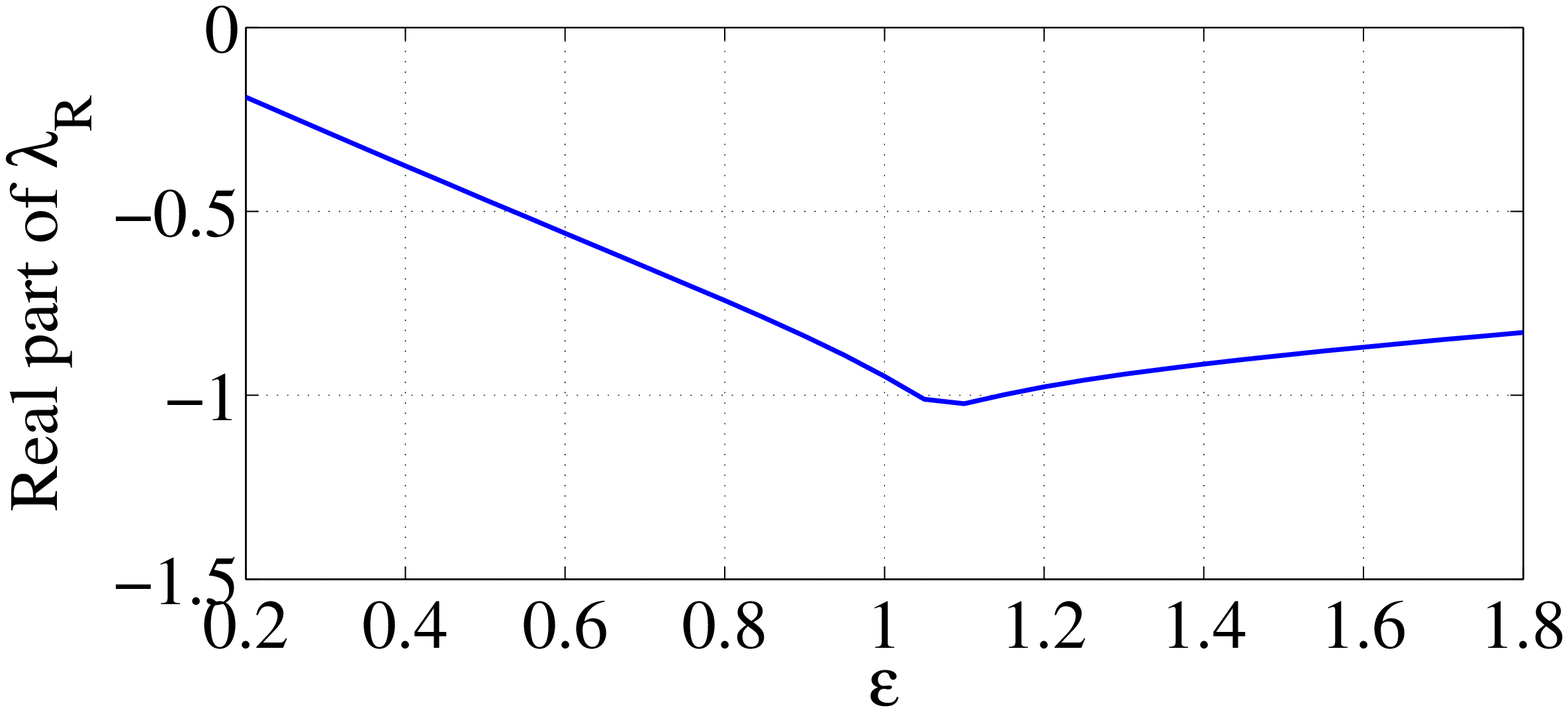}
\caption{Real part of the rightmost non-null eigenvalue of matrix $M(\varepsilon)$ with respect to $\varepsilon$.}\label{optimalepsilon}
\end{figure}

\begin{figure}
\centering
\includegraphics[width=0.9\linewidth]{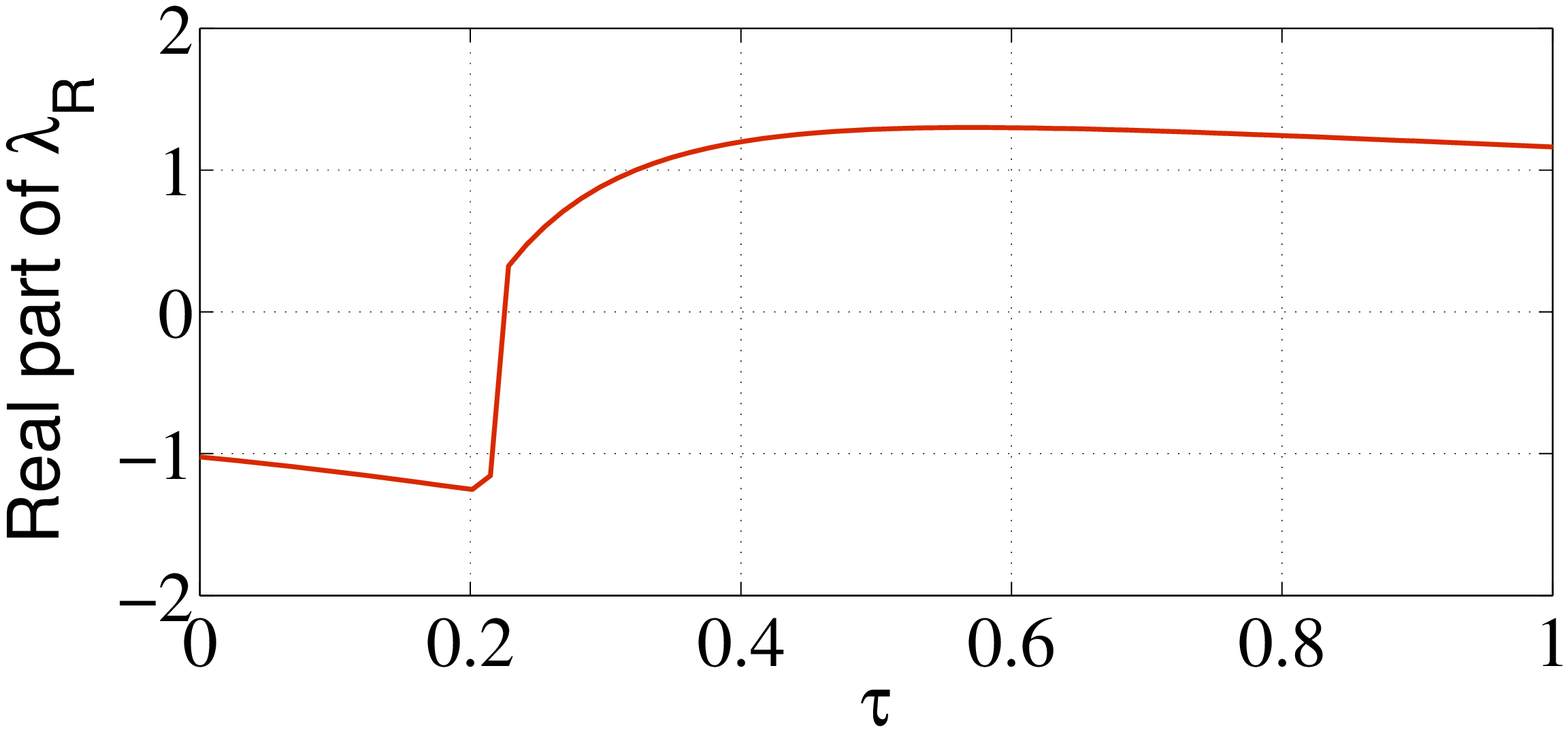}
\caption{Real part of rightmost non-null root of eq.~\eqref{charpoly} with respect to $\tau$, for $\varepsilon=1.1$.}\label{fig:optimumeps}
\end{figure}

\begin{figure}
\centering
\includegraphics[width=0.9\linewidth]{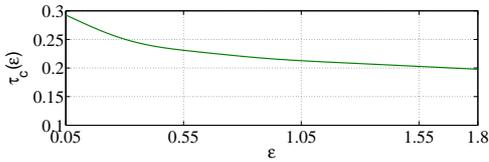}
\caption{The value of $\tau_c(\varepsilon)$ with respect to $\varepsilon$.}\label{crossing}
\end{figure}

\section{Numerical example and simulations}\label{section:example}

In this section we consider a numerical example to corroborate the theoretical results presented in the previous section.

We consider the network of $6$ agents whose topology is shown in Fig.~\ref{topology}. Such a network is encoded by the adjacency matrix
\begin{equation}\label{adjacency_mat}
A=\left[\begin{array}{cccccc}
0 & 1 & 0 & 1 & 0 & 1  \\
0 & 0 & 1 & 0 & 1 & 0  \\
0 & 1 & 0 & 1 & 0 & 1  \\
1 & 0 & 0 & 0 & 1 & 0 \\
0 & 1 & 0 & 0 & 0 & 1 \\
1 & 0 & 1 & 0 & 0 & 0
\end{array} \right]
\end{equation}

The in and out-Laplacian matrices are, respectively
$$\mathcal{L}_{in}=\left[\begin{array}{cccccc}
3  & -1 &  0 & -1 & 0 & -1  \\
0  &  2 & -1 & 0 & -1 & 0  \\
0  & -1 &  3 & -1 & 0 & -1  \\
-1 &  0 &  0 & 2 & -1 & 0 \\
0  & -1 &  0 & 0 & 2 & -1 \\
-1 &  0 & -1 & 0 & 0 & 2
\end{array} \right]
$$
and
$$\mathcal{L}_{out}=\left[\begin{array}{cccccc}
2  & -1 &  0 & -1 & 0 & -1  \\
0  &  3 & -1 & 0 & -1 & 0  \\
0  & -1 &  2 & -1 & 0 & -1  \\
-1 &  0 &  0 & 2 & -1 & 0 \\
0  & -1 &  0 & 0 & 2 & -1 \\
-1 &  0 & -1 & 0 & 0 & 3
\end{array} \right].
$$

 Fig.~\ref{fig:examplex} shows the evolution of system~\eqref{withdelay} when $\varepsilon=1.3$ and $\tau=0.18$. Initial conditions $x(0)$ are chosen uniformly at random while initial conditions $z(0)=\mathbf{0}_n$. Fig.~\ref{fig:examplex} shows how consensus on the average of the initial state $x(0)$ is achieved. Fig.~\eqref{fig:examplez} presents the evolution of the storage variables $z(t)$. All storage variables are initially set to zero and then vary during the dynamical evolution of the system so that the quantity $\mathbf{1}_n^Tx(t)+\mathbf{1}_n^Tz(t)$ remains constant.

\begin{figure}
\centering
\includegraphics[width=0.9\linewidth]{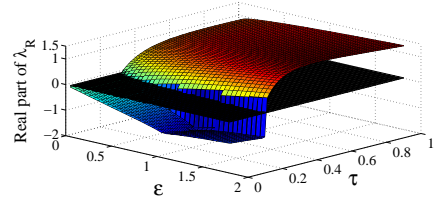}
\caption{Value of the real part of the rightmost non-null root $\lambda_R$ of eq.~\eqref{charpoly} versus increasing $\varepsilon$ and time delay $\tau$.  }\label{fig:tau01eps001}
\end{figure}

We now present the results of a series of numerical simulations whose aim is that of showing how the consensus achievement is related to parameters $\varepsilon$ and $\tau$. In particular, Fig.~\ref{optimalepsilon} shows how the rightmost non-null eigenvalue $\lambda_R$ of matrix $M(\varepsilon)$ varies for $\varepsilon \in \left[0.2,1.8\right]$. Fig.~\ref{optimalepsilon} shows that there exists an optimal value at $\varepsilon=1.1$ for which matrix $M(\varepsilon)$ in the given example has the smallest rightmost non-null eigenvalue.

In Fig.~\ref{fig:optimumeps} we show how the rightmost non-null root of eq.~\eqref{charpoly} varies for increasing values of the time-delay $\tau$ when $\varepsilon=1.1$. Fig.~\ref{fig:optimumeps} shows that despite the time-delay can make the system unstable, it can also improve the convergence speed to average consensus. For this example the optimal value of the time-delay is $\tau=0.19$.

Fig.~\ref{crossing} shows the values of $\tau_c$ in eq.~\eqref{tau_c} for which eq.~\eqref{charpoly} has roots in the imaginary axis, i.e., it shows the maximum time delay sustainable by system in eq.~\eqref{withdelay} for the considered network topology in Fig.~\ref{topology}.

Finally, in Fig.~\ref{fig:tau01eps001} we show a plot of the real part of the rightmost non-null eigenvalue of eq.~\eqref{charpoly} for $\varepsilon\in \left(0,2\right]$ and $\tau\in \left[0,1\right]$. Fig.~\ref{fig:tau01eps001} shows how the convergence properties are affected by parameters $\varepsilon$ and $\tau$: there exists an optimal value at $\varepsilon=1.1$ and $\tau=0.19$ for which $\lambda_R$ is the most negative and there exists a connected region of the plane defined by $\varepsilon,\tau$ where $\lambda_R$ has strictly negative real part.

The rightmost non-null root of eq.~\eqref{charpoly} for a given set of $(\varepsilon, \tau)$ is computed using the spectral method with the heuristic presented by \cite{wu2012reliably}. 


\section{Conclusions}\label{section:conclusion}
In this preliminary paper a continuous time version of a consensus on the average protocol for arbitrary strongly connected directed graphs is proposed and its convergence properties with respect to time delays in the local state update are characterized. The convergence properties of this algorithm depend upon a tuning parameter that can be made arbitrary small to prove stability of the networked system. Simulations have been presented to corroborate the theoretical results and show that the existence of a small time delay can actually improve the algorithm performance.

Future work will include an extension of the mathematical characterization of the proposed algorithm to consider possibly heterogeneous or time-varying delays.

\bibliographystyle{plain}
\bibliography{biblio}

\end{document}